\documentclass[12pt,twoside]{amsart}
\usepackage{amssymb}
\usepackage{amscd}
\usepackage{amsmath}
\usepackage[all]{xy}

\title{Semipositivity theorems for moduli problems}
\author{Osamu Fujino}
\date{2018/2/11, version 2.01}
\address{Department of Mathematics, Graduate School of Science,
Osaka University, Toyonaka, Osaka 560-0043, Japan}
\email{fujino@math.sci.osaka-u.ac.jp}

\makeatletter
    
    \@addtoreset{equation}{section}
\makeatother
\newcommand{\Sing}[0]{\operatorname{Sing}}
\newcommand{\Gr}[0]{\operatorname{Gr}}
\newcommand{\Supp}[0]{\operatorname{Supp}}
\newcommand{\Aut}[0]{\operatorname{Aut}}
\newcommand{\Spec}[0]{\operatorname{Spec}}
\newcommand{\Exc}[0]{\operatorname{Exc}}
\newcommand{\codim}[0]{\operatorname{codim}}

\newtheorem{thm}{Theorem}[section]
\newtheorem{lem}[thm]{Lemma}
\newtheorem{cor}[thm]{Corollary}

\theoremstyle{definition}

\newtheorem{defn}[thm]{Definition}
\newtheorem{rem}[thm]{Remark}
\newtheorem*{ack}{Acknowledgments}

\newtheorem{say}[thm]{}
\newtheorem{step}{Step}
\begin{document}
\bibliographystyle{amsalpha+}

\maketitle

\begin{abstract}
We prove some semipositivity theorems for singular varieties coming from
graded polarizable admissible variations of mixed Hodge structure.
As an application, we obtain that the moduli functor of stable varieties 
is semipositive in the sense of Koll\'ar. This completes Koll\'ar's 
projectivity criterion for the moduli spaces of higher-dimensional 
stable varieties. 
\end{abstract}

\tableofcontents
\section{Introduction}\label{ff-sec1}

The main 
purpose of this paper is to give a proof of the following ``folklore statement"
(see, for example, \cite{alexeev2}, \cite{karu}, 
\cite{kovacs}, \cite{abramovich-karu}, \cite{kollar2}, and \cite{fujino-sugaku}) 
based on \cite{fujino-fujisawa} (see also \cite{ffs}). 
In general, the (quasi-)projectivity 
of some moduli space is a subtle problem and is
harder than it looks (see, 
for example, \cite{kollar-ann}, and \cite{viehweg3}). 
We note that the coarse moduli space of stable varieties was first 
constructed in the category of 
algebraic spaces. 

\begin{thm}[Projectivity of moduli spaces of stable varieties]\label{ff-thm1.1}
Every closed complete subspace of the 
coarse moduli space of stable varieties is {\em{projective}}.
\end{thm}

To the best knowledge of the author, Theorem \ref{ff-thm1.1} 
is new for stable $n$-folds with $n\geq 3$ (see 
the comments in \ref{ff-say1.8} below). 
Note that a {\em{stable $n$-fold}} is an $n$-dimensional 
projective semi-log-canonical variety with ample canonical 
divisor and is 
called an {\em{$n$-dimensional semi-log-canonical 
model}} in \cite{kollar2}. For the details of semi-log-canonical varieties, 
see, for example, \cite{fujino} and \cite{kollar-book}. 
Theorem \ref{ff-thm1.1} is a direct 
consequence of Theorem \ref{ff-thm1.2} below
by Koll\'ar's projectivity criterion (see \cite[Sections 2 and 3]{kollar}).

\begin{thm}[Semipositivity of $\mathcal M^{stable}$]\label{ff-thm1.2}
Let $\mathcal M^{stable}$ be the moduli functor of stable varieties. 
Then $\mathcal M^{stable}$ 
is semipositive in the sense of Koll\'ar. 
\end{thm}

For the reader's convenience, let us recall the 
definition of the semipositivity of $\mathcal 
M^{stable}$, which is a special case of \cite[2.4.~Definition]{kollar}. 

\begin{defn}[{see \cite[2.4.~Definition]{kollar}}]\label{ff-def1.3} 
The moduli functor $\mathcal M^{stable}$ of stable varieties is said to be 
{\em{semipositive}} 
({\em{in the sense of Koll\'ar}}) if the following condition holds: 

There is a fixed positive integer 
$m_0$ such that if $C$ is a smooth projective curve and $(f:X\to C)\in 
\mathcal M^{stable}(C)$, 
then $f_*\omega^{[mm_0]}_{X/C}$ is a nef locally 
free sheaf on $C$ for every positive integer $m$. 
\end{defn}

As the culmination of the works of several authors 
(see, for example, \cite{alexeev1}, 
\cite{alexeev2}, \cite{alexeev-mori}, 
\cite{hmx}, 
\cite{hacon-xu}, \cite{keel-mori}, 
\cite{kollar-hull}, \cite{kollar-s}, and \cite{ksb}), we have: 
 
\begin{cor}\label{ff-cor1.4} 
The moduli functor $\mathcal M^{stable}_H$ of 
stable varieties with Hilbert function $H$ 
is coarsely represented by a {\em{projective}} 
algebraic scheme. 
\end{cor}

As an easy consequence of 
Corollary \ref{ff-cor1.4}, we obtain:

\begin{cor}[{see \cite[Theorem 1.11]{viehweg2}}]\label{ff-cor1.5} 
The moduli functor $\mathcal M_H$ of canonically polarized smooth projective 
varieties with Hilbert function $H$ is coarsely represented by 
a {\em{quasi-projective}} algebraic scheme. 
\end{cor}

More generally, we have: 

\begin{cor}\label{ff-cor1.6} 
The moduli functor $\mathcal M^{can}_H$ of 
canonically polarized normal projective varieties 
having only canonical singularities 
with Hilbert function $H$ is coarsely represented by a {\em{quasi-projective}} 
algebraic scheme. 
\end{cor}

Theorem \ref{ff-thm1.2} follows almost directly from the definition of the semipositivity of 
$\mathcal M^{stable}$ in Definition \ref{ff-def1.3} 
(see \cite[2.4.~Definition]{kollar}) and the following semipositivity theorem. 

\begin{thm}[Semipositivity theorem I]
\label{ff-thm1.7}
Let $X$ be an equidimensional variety 
which satisfies Serre's $S_2$ condition and is normal crossing in 
codimension one.
Let $f:X\to C$ be a projective
surjective morphism onto a smooth projective curve $C$ such that
every irreducible component of
$X$ is dominant onto $C$.
Assume that there exists a non-empty Zariski open set $U$ of $C$ such that
$f^{-1}(U)$ has only semi-log-canonical singularities.
Then $f_*\omega_{X/C}$ is nef. 

Assume further that $\omega_{X/C}^{[k]}$ is locally free and $f$-generated for some positive
integer $k$.
Then $f_*\omega_{X/C}^{[m]}$ is nef for every $m\geq 1$.
\end{thm}

\begin{say}[Comments]\label{ff-say1.8}  
By the recent developments of the minimal model program, 
Theorem \ref{ff-thm1.7} seems to be 
a reasonable formulation of \cite[4.12.~Theorem]{kollar} 
for higher-dimensional singular varieties.
Koll\'ar has pointed out that
the assumption that the fibers are surfaces was
inadvertently omitted from its statement.
He is really claiming \cite[4.12.~Theorem]{kollar} for $f:Z\to C$ with $\dim Z=3$ 
(see \cite[1.~Introduction]{kollar}).
Therefore,
Theorem \ref{ff-thm1.7} is new when $\dim X\geq 4$.
Likewise, Theorems \ref{ff-thm1.1} and \ref{ff-thm1.2} 
are new when the dimension of the
stable varieties are greater than or equal to three. 
We feel that the arguments in \cite[4.14]{kollar}
only work when the fibers are surfaces.
In other words, we needed some new ideas and techniques
to prove Theorem \ref{ff-thm1.7}.
Our arguments heavily depend on the recent 
advances on the semipositivity theorems of
Hodge bundles (\cite{fujino-fujisawa} and 
\cite{ffs}) and some ideas in \cite{fujino}. 
\end{say} 

For the general theory of Koll\'ar's projectivity criterion,
see \cite[Sections 2 and 3]{kollar} and \cite[Theorem 4.34]{viehweg2}.
We do not discuss the technical details of the construction of moduli spaces 
of stable varieties in this paper. We mainly treat various semipositivity theorems. 
Note that the projectivity criterion discussed here is independent of 
the existence problem of moduli spaces. 
We recommend the reader to see \cite[Section 2]{kollar} 
and \cite[Sections 4 and 5]{kollar2} for Koll\'ar's program for
constructing
moduli spaces of stable varieties (see also \cite{kollar-ongoing}).
Our paper is related to the topic in \cite[5.5 (Projectivity)]{kollar2}.

In this paper, we prove Theorem \ref{ff-thm1.7}
in the framework of \cite{fujino} and \cite{fujino-fujisawa}, 
although we do not use the arguments in \cite{fujino} explicitly.
Note that \cite{fujino} 
and \cite{fujino-fujisawa} heavily depend on the theory of {\em{mixed}} 
Hodge structures
on cohomology with compact support.
A key ingredient of this paper is the following semipositivity theorem,
which is essentially contained in \cite{fujino-fujisawa} (see 
also \cite{ffs}).
It is a generalization of Fujita's semipositivity theorem (see 
\cite[(0.6) Main Theorem]{fujita}). 
We note that a Hodge theoretic approach to 
the original Fujita semipositivity theorem 
was introduced by Zucker (see \cite{zucker}). 

\begin{thm}[{Basic semipositivity theorem, 
see \cite[Section 7]{fujino-fujisawa}}]\label{ff-thm1.9}
Let $(X, D)$ be a simple normal crossing pair such that $D$ is reduced.
Let $f:X\to C$ be a projective surjective morphism onto a smooth projective
curve $C$.
Assume that every stratum of $X$ is dominant onto $C$.
Then $f_*\omega_{X/C}(D)$ is nef.
\end{thm}

Although we do not know what is the best 
formulation of the semipositivity theorem
for moduli problems, we think Theorem \ref{ff-thm1.9} will 
be one of the most fundamental
results for application to  Koll\'ar's projectivity criterion for moduli spaces.
We can prove Theorem \ref{ff-thm1.7} by using Theorem \ref{ff-thm1.9}.
By the same proof as that of Theorem \ref{ff-thm1.7},
we obtain a generalization of Theorem \ref{ff-thm1.7},
which implies both Theorems \ref{ff-thm1.7} and \ref{ff-thm1.9}.

\begin{thm}[Semipositivity theorem II]\label{ff-thm1.10}
Let $X$ be an equidimensional variety which 
satisfies Serre's $S_2$ condition and
is normal crossing in codimension one.
Let $f:X\to C$ be a projective surjective morphism onto 
a smooth projective curve $C$ such that
every irreducible component of $X$ is dominant onto $C$.
Let $D$ be a reduced Weil divisor on $X$ such that
no irreducible component of $D$ is contained in the singular locus of $X$.
Assume that there exists a non-empty 
Zariski open set $U$ of $C$ such that $(f^{-1}(U), D|_{f^{-1}(U)})$ is
a semi-log-canonical pair.
Then $f_*\omega_{X/C}(D)$ is nef.

We further assume that $\mathcal O_X(k(K_X+D))$ is locally free and $f$-generated
for some positive integer $k$. Then
$f_*\mathcal O_X(m(K_{X/C}+D))$ is nef for every $m\geq 1$.
\end{thm}

By combining Theorem \ref{ff-thm1.10} with 
Viehweg's covering trick, we obtain the 
following theorem:~Theorem \ref{ff-thm1.11}, 
which is an answer to the question in \cite[5.6]{alexeev2}. 
Although we do not discuss the moduli spaces of {\em{stable pairs}} here,
Theorems \ref{ff-thm1.10} and \ref{ff-thm1.11} 
play important roles in the proof of the projectivity of the
moduli spaces of stable pairs (see \cite{alexeev2}, \cite{fp}, 
\cite{hassett}, \cite{kp}, and \ref{ff-say4.3} below).

\begin{thm}[Semipositivity theorem III]\label{ff-thm1.11}
Let $X$ be an equidimensional variety 
which satisfies Serre's $S_2$ condition and
is normal crossing in codimension one.
Let $f:X\to C$ be a projective surjective 
morphism onto a smooth projective curve $C$ such that
every irreducible component of $X$ is dominant onto $C$.
Let $\Delta$ be an effective $\mathbb Q$-Weil divisor on $X$ such that
no irreducible component of the support 
of $\Delta$ is contained in the singular locus of $X$.
Assume that there exists a non-empty Zariski 
open set $U$ of $C$ such that $(f^{-1}(U), \Delta|_{f^{-1}(U)})$ is
a semi-log-canonical pair. 
We further assume that $\mathcal O_X(k(K_X+\Delta))$ is locally free and 
$f$-generated for some positive integer $k$. 
Then $f_*\mathcal O_X(k(K_{X/C}+\Delta))$ is nef. 
Therefore, $f_*\mathcal O_X(kl(K_{X/C}+\Delta))$ is nef for every $l\geq 1$. 
\end{thm}

Recently, Kov\'acs and Patakfalvi generalized 
Theorem \ref{ff-thm1.1} for stable 
pairs in \cite{kp}. 
Note that one of the main ingredients of \cite{kp} is Theorem \ref{ff-thm1.11}. 
For the details, we recommend 
the reader to see \cite{kp}. 
Theorem \ref{ff-thm1.11} is also a key result for the proof of the ampleness of 
the CM line bundle on the moduli space of canonically 
polarized varieties in \cite{px}. 
In any case, we expect our semipositivity 
theorems established in \cite{fujino-fujisawa} 
to play an important role in the study of higher-dimensional 
complex algebraic varieties. 

\begin{rem}\label{ff-rem1.12}
In this paper, we do not use algebraic spaces for the proof of the semipositivity theorems.
We only treat {\em{projective}} varieties.
Note that Theorem \ref{ff-thm1.9} follows from the theory of variations of mixed Hodge
structure. The variations of mixed Hodge structure discussed 
in \cite{fujino-fujisawa} (see also \cite{ffs}) 
are {\em{graded polarizable}} and {\em{admissible}}.
Therefore, we can not directly apply the results in \cite{fujino-fujisawa} to
the variations of (mixed) Hodge structure arising from families of
algebraic spaces. We need
some polarization to obtain various semipositivity theorems in our framework.
We also note that the admissibility assures us the existence of canonical extensions of Hodge bundles, which
does not always hold for abstract graded polarizable variations of
mixed Hodge structure (see \cite[Example 1.5]{fujino-fujisawa}).
\end{rem}

We do not use the Fujita--Zucker--Kawamata semipositivity
theorem coming from the theory of polarized variations of Hodge structure 
(see \cite{fujino-zucker65}).

\begin{rem}\label{ff-rem1.13}
(1) As explained in \cite{kollar} and \cite{kollar-two}, it is difficult to 
directly check the {\em{quasi-projectivity}} of non-complete singular spaces. 
This is because there is no good ampleness criterion for non-complete spaces. 
In this paper, we adopt Koll\'ar's framework in \cite[Sections 2 and 3]{kollar}, where 
we use the Nakai--Moishezon criterion to check 
the projectivity of complete algebraic spaces. 
Note that Viehweg discusses the {\em{quasi-projectivity}} of {\em{non-complete}} 
moduli spaces (see \cite{viehweg2} and \cite{viehweg3}). On the other hand, 
Koll\'ar and we prove the {\em{projectivity}} 
of {\em{complete}} moduli spaces (see \cite{kollar}). 

(2) In general, we have to formulate and prove semipositivity theorems 
for non-normal (reducible) 
varieties $X$ even if we are mainly interested in the moduli spaces of smooth 
(or normal) stable varieties. 
Let $\mathcal M_g$ (resp.~$\overline{\mathcal M}_g$) 
be the moduli functor of smooth projective 
curves (resp.~stable curves) with $g\geq 2$. 
We consider $(f:X\to C) \in \overline{\mathcal M}_g(C)$ 
such that 
$C$ is contained in $\overline{M}_g\setminus M_g$, where 
$M_g$ (resp.~$\overline{M}_g$) is the coarse moduli 
space of $\mathcal M_g$ (resp.~$\overline{\mathcal M}_g$). 
Then a general fiber of $f:X\to C$ may be non-normal and reducible. 
We can prove the projectivity of $\overline{M}_g$ by Koll\'ar's 
projectivity criterion. However, we can not directly prove 
the quasi-projectivity of $M_g$. 

(3) Although we repeatedly use 
Viehweg's covering arguments, we do not use the notion of 
{\em{weak positivity}}, which was introduced by Viehweg and plays 
crucial roles in his works (see \cite{viehweg1}, \cite{viehweg2}, and 
\cite{viehweg3}). 
We just treat the {\em{semipositivity}} on smooth projective curves (see \cite{kollar}). 

(4) From the Hodge theoretic viewpoint, 
our approach is based on the theory of {\em{mixed}} 
Hodge structures (see \cite{fujino-fujisawa} and \cite{ffs}). 
The arguments in \cite{viehweg2}, \cite{kollar}, 
and \cite{viehweg3} use only {\em{pure}} 
Hodge structures. It is one of the main differences 
between our approach and the others. 

(5) In this paper, we use the theory of variations of mixed Hodge 
structure only in the proof of Theorem \ref{ff-thm1.9}. 
Moreover, for the proof of Theorem \ref{ff-thm1.9}, 
we only need the theory of variations of mixed Hodge structure 
in the case where the base space is a curve. 
If we assume that the base space is a curve, 
then the theory described in \cite{fujino-fujisawa} becomes 
much simpler than the general case. 
\end{rem}

We summarize the contents of this paper. 
In Section \ref{ff-sec2}, 
we collect some basic definitions. 
In Section \ref{ff-sec3}, we quickly 
review the moduli functor $\mathcal M^{stable}$ of 
stable varieties and its coarse moduli space. 
Section \ref{ff-sec4} is the main part of this paper, where 
we prove the theorems in Section \ref{ff-sec1}. 
Our proofs depend on \cite{fujino-fujisawa}, some ideas in 
\cite{fujino}, 
and Viehweg's covering arguments. 
In Section \ref{ff-sec5}, we prove the corollaries in Section \ref{ff-sec1}. 

\begin{ack}
The author was partially supported by JSPS KAKENHI Grant Numbers 
JP16H03925, JP16H06337. 
When he wrote the original version of this paper in 2012, 
he was partially supported 
by Grant-in-Aid for Young Scientists (A) 24684002 from JSPS.
He thanks Professor J\'anos Koll\'ar for answering 
his questions and giving him many useful comments 
and Professor Steven Zucker for giving him useful 
comments and advice.
He also thanks Professors Varely Alexeev, Taro 
Fujisawa, S\'andor Kov\'acs, and referees for comments. 
Finally, he thanks Professor Shigefumi Mori for warm encouragement. 
\end{ack}

We will work over $\mathbb C$, the complex 
number field, throughout this paper.
Note that, by the Lefschetz principle,
all the results in this paper hold over any 
algebraically closed field $k$ of characteristic zero.
We will freely use the notation and terminology in \cite{fujino-fujisawa} and \cite{fujino}. 
For the standard 
notations and conventions of the 
log minimal model program, see \cite{fujino-fund} and \cite{fujino-foundation}. 

\section{Preliminaries}\label{ff-sec2}
Let us recall the definition of {\em{nef locally free sheaves}}.
For the details, see, for example, \cite[Section 2]{viehweg2}.

\begin{defn}[Nef locally free sheaves]\label{ff-def2.1}
A locally free sheaf of finite rank $\mathcal E$ on a complete
variety $X$ is {\em{nef}} if the following equivalent conditions are satisfied:
\begin{itemize}
\item[(i)] $\mathcal E=0$ or $\mathcal O_{\mathbb P_X(\mathcal E)}(1)$ is nef
on $\mathbb P_X(\mathcal E)$.
\item[(ii)] For every map from a smooth projective curve $f:C\to X$,
every quotient line bundle of $f^*\mathcal E$ has non-negative degree.
\end{itemize}
A nef locally free sheaf was originally called a ({\em{numerically}}) 
{\em{semipositive}} sheaf in the literature. 
\end{defn}

In this paper, we only discuss various semipositivity theorems
for locally free sheaves on a smooth projective curve. 
The following well-known lemma is very useful. 
We omit the proof of Lemma \ref{ff-lem2.2} because it is an easy exercise.

\begin{lem}\label{ff-lem2.2}
Let $C$ be a smooth projective curve and let $\mathcal E_i$ be
a locally free sheaf on $C$ for $i=1, 2$.
Assume that $\mathcal E_1\subset \mathcal E_2$, $\mathcal E_1$ is
nef, and that $\mathcal E_1$ coincides with $\mathcal E_2$ over some
non-empty Zariski open set of $C$.
Then $\mathcal E_2$ is nef.
\end{lem}

The following lemma is more or less known to the 
experts. However, we can not find it explicitly in the literature. 
So we describe it here for the reader's convenience. 

\begin{lem}\label{ff-lem2.3}
Let $\mathcal E$ be a locally free sheaf of finite rank on a smooth 
projective 
irreducible curve $C$. 
Let $\Sigma$ be a fixed Zariski closed set of $C$ with $\Sigma\subsetneq C$. 
Assume that there exists some positive integer $\mu$ such that, for 
every finite surjective 
morphism $\pi:C'\to C$ from a smooth projective 
irreducible curve $C'$ which is 
\'etale over some open neighborhood of $\Sigma$ and 
for every ample line bundle $\mathcal H'$ on $C'$, 
$\pi^*\mathcal E\otimes (\mathcal H')^\mu$ is nef on $C'$. 
Then $\mathcal E$ is nef.   
\end{lem}
\begin{proof}
We take an ample line bundle $\mathcal H$ on $C$. 
For any positive integer $\alpha$, we can construct a finite 
covering $\pi:C'\to C$ from a smooth projective 
irreducible curve $C'$ such that 
$\pi^*\mathcal H=(\mathcal H')^{1+2\alpha \mu}$ for 
some ample line bundle $\mathcal H'$ on $C'$. 
We can make $\pi$ \'etale over some open neighborhood 
of $\Sigma$ and assume that 
$\pi$ is Galois (see \cite[Theorem 1-1-1]{kmm}). 
We note that the trace map splits the natural inclusion 
$\mathcal O_C\to \pi_*\mathcal O_{C'}$. 
Since $\pi^*\mathcal E\otimes (\mathcal H')^\mu$ is nef, 
there is some positive integer $\beta$ such that 
$$
S^{(2\alpha)\beta}(\pi^*\mathcal E\otimes (\mathcal H')^\mu)\otimes 
(\mathcal H')^\beta=\pi^*(S^{2\alpha\beta}(\mathcal E)\otimes 
\mathcal H^\beta)
$$ 
is generated by its global sections (see \cite[Proposition 2.9]{viehweg2}). 
Therefore, we have a surjective 
morphism 
$$
\bigoplus _{\mathrm{finite}} \mathcal O_{C'}\to \pi^*(S^{2\alpha\beta}
(\mathcal E)\otimes \mathcal H^\beta). 
$$
Thus the induced morphism 
$$
\left(\bigoplus _{\mathrm{finite}}\pi_*\mathcal O_{C'} 
\right)\otimes \mathcal H^\beta\to 
S^{2\alpha\beta}(\mathcal E)\otimes \mathcal H^{2\beta}\otimes 
\pi_*\mathcal O_{C'}\to S^{2\alpha\beta}(\mathcal E)\otimes 
\mathcal H^{2\beta}
$$ 
is surjective. By replacing $\beta$ by some multiple, we may assume 
that $\pi_*\mathcal O_{C'}\otimes \mathcal H^\beta$ is generated 
by its global sections. 
In this case, $S^{2\alpha\beta}(\mathcal E)\otimes \mathcal H^{2\beta}$ is generated 
by its global sections. 
This implies that $\mathcal E$ is nef (see \cite[Proposition 2.9]{viehweg2}). 
\end{proof}

We need the notion of {\em{simple normal crossing pairs}}
for Theorem \ref{ff-thm1.9}. Note that a simple normal crossing pair is sometimes
called a {\em{semi-snc}} pair in the literature (see \cite[Definition 1.1]{bierstone}).

\begin{defn}[Simple normal crossing pairs]\label{ff-def2.4}
We say that the pair $(X, D)$ is {\em{simple normal crossing}} at
a point $a\in X$ if $X$ has a 
Zariski open neighborhood $U$ of $a$ that can be embedded in a smooth
variety
$Y$,
where $Y$ has regular system of parameters 
$(x_1, \cdots, x_p, y_1, \cdots, y_r)$ at
$a=0$ in which $U$ is defined by a monomial equation
$$
x_1\cdots x_p=0
$$
and $$
D=\sum _{i=1}^r \alpha_i(y_i=0)|_U, \quad  \alpha_i\in \mathbb R.
$$
We say that $(X, D)$ is a {\em{simple normal crossing pair}} if it is 
simple normal crossing at every point of $X$. 
We sometimes 
say that $D$ is a {\em{simple normal crossing divisor}} on $X$ 
if $(X, D)$ is a simple normal crossing pair and $D$ is reduced. 
If $(X, 0)$ is a simple normal crossing pair, 
then we simply say that $X$ is a {\em{simple normal crossing variety}}. 
Let $X$ be a simple normal crossing variety and let $X=\sum _{i\in I} X_i$ 
be the irreducible decomposition. 
A {\em{stratum}} of $X$ is an irreducible component of 
$X_{i_1}\cap 
\cdots \cap X_{i_k}$ for some $\{ i_1, \cdots, i_k\}\subset I$. 
\end{defn}

\begin{defn}[Stratum]\label{ff-def2.5}
Let $(X, D)$ be a simple normal crossing pair such that 
$D$ is reduced. 
Let $\nu:X^\nu\to X$ be the normalization. 
We put $K_{X^\nu}+\Theta=\nu^*(K_X+D)$, 
that is, $\Theta$ is the sum of the inverse images of $D$ and 
the singular locus of $X$. 
A {\em{stratum}} of $(X, D)$ is an irreducible component 
of $X$ or the $\nu$-image of a log canonical 
center of $(X^\nu, \Theta)$. 
This definition is compatible with Definition \ref{ff-def2.4}. 
\end{defn}

For the reader's convenience, we recall the notion of {\em{semi-log-canonical pairs}}. 

\begin{defn}[Semi-log-canonical pairs]\label{ff-def2.6} 
Let $X$ be an equidimensional algebraic variety which satisfies Serre's 
$S_2$ condition and is normal crossing in codimension one. 
Let $\Delta$ be an effective $\mathbb R$-divisor 
on $X$ such that no irreducible component of $\Supp \Delta$ is contained 
in the singular locus of $X$. 
The pair $(X, \Delta)$ is called a {\em{semi-log-canonical pair}} 
(an {\em{slc pair}}, for short) if 
\begin{itemize}
\item[(1)] $K_X+\Delta$ is $\mathbb R$-Cartier, and 
\item[(2)] $(X^\nu, \Theta)$ is log canonical, where 
$\nu:X^\nu\to X$ is the normalization and 
$K_{X^\nu}+\Theta=\nu^*(K_X+\Delta)$, that is, 
$\Theta$ is the sum of the inverse images of $\Delta$ and 
the conductor of $X$. 
\end{itemize}
If $(X, 0)$ is a semi-log-canonical pair, then we simply say that 
$X$ is a {\em{semi-log-canonical variety}} or $X$ has only {\em{semi-log-canonical 
singularities}}. 
\end{defn}

For the details of semi-log-canonical pairs and 
the basic notations, see \cite{fujino} and \cite{kollar-book}. 

\begin{say}[$\mathbb Q$-divisors]\label{ff-say2.7} 
Let $D$ be a $\mathbb Q$-divisor 
on an equidimensional variety $X$, that is, 
$D$ is a finite formal $\mathbb Q$-linear 
combination 
$$D=\sum _i d_i D_i$$ of irreducible 
reduced subschemes $D_i$ of codimension one such 
that $D_i\ne D_j$ for $i\ne j$. 
We define the {\em{round-up}} 
$\lceil D\rceil =\sum _i \lceil d_i \rceil D_i$ (resp.~{\em{round-down}} 
$\lfloor D\rfloor =\sum _i \lfloor d_i \rfloor D_i$), where 
every real number $x$, $\lceil x\rceil$ (resp.~$\lfloor x\rfloor$) is the integer 
defined by $x\leq  \lceil x \rceil <x+1$ 
(resp.~$x-1<\lfloor x\rfloor \leq x$). 
The {\em{fractional part}} $\{D\}$ of $D$ denotes $D-\lfloor D\rfloor$. 
We set 
$$D^{<0}=\sum _{d_i<0}d_i D_i,  \quad 
D^{>0}=\sum _{d_i>0}d_i D_i, \quad \text{and}\quad D^{=1}=\sum _{d_i=1}D_i.$$ 
\end{say}

\begin{say}[Demi-normal variety]\label{ff-say2.8} 
Let $X$ be an equidimensional variety which 
satisfies Serre's $S_2$ condition and is normal crossing 
in codimension one. 
Then $X$ is sometimes said to be {\em{demi-normal}} 
(see \cite[Definition 5.1]{kollar-book}). 

Let $\pi:Y\to X$ be a finite surjective morphism 
between demi-normal varieties and let $D$ be a $\mathbb Q$-divisor 
on $X$ such that 
no irreducible component of $D$ is contained in the singular 
locus of $X$. Then there is a Zariski open set $U$ of $X$ such that 
$\codim _X(X\setminus U)\geq 2$ and that 
$D|_U$ is $\mathbb Q$-Cartier on $U$. 
Assume that $\pi$ is \'etale over 
the generic point of any irreducible component 
of $\Supp D$. 
Then we have a well-defined 
$\mathbb Q$-Cartier $\mathbb Q$-divisor $\pi^*(D|_U)$ on 
$\pi^{-1}(U)$. In this situation, 
$\pi^{-1}D$ denotes the $\mathbb Q$-divisor on $Y$ that is the closure of 
$\pi^*(D|_U)$. By construction, no irreducible component of $\pi^{-1}D$ is 
contained in the singular locus of $Y$. 
\end{say}

We recall the definition of $\omega_{X/C}^{[m]}$.

\begin{defn}\label{ff-def2.9}
In Theorem \ref{ff-thm1.7},
$\omega_{X/C}^{[m]}$ is the {\em{$m$-th
reflexive
power}} of $\omega_{X/C}$.
It is the double dual of the $m$-th tensor power
of $\omega_{X/C}$:
$$
\omega_{X/C}^{[m]}:=(\omega_{X/C}^{\otimes m})^{**}.
$$
\end{defn}

For the details of divisors and divisorial sheaves on demi-normal varieties, 
see \cite[Section 5.1]{kollar-book}. 

\begin{say}[Lemmas on resolution of singularities for reducible 
varieties]\label{ff-say2.10}
We prepare the following lemmas on resolution of singularities. 
We will use them in Section \ref{ff-sec4}. 

\begin{lem}\label{ff-lem2.11}
Let $(Y, \Delta)$ be a simple normal crossing 
pair. 
Let $D$ be a Weil {\em{(}}resp.~$\mathbb Q$-Weil or $\mathbb R$-Weil{\em{)}} 
divisor on $Y$ such that 
$\Supp D\subset \Supp \Delta$. 
Then there exists a sequence of blow-ups 
$$
Y=Y_0\overset{\pi_1}{\longleftarrow} 
Y_1\overset{\pi_2}{\longleftarrow} Y_2
\overset{\pi_3}{\longleftarrow} \cdots \overset{\pi_k}{\longleftarrow} Y_k
$$
with the following properties. 
\begin{itemize}
\item[(1)] $(Y_0, \Delta_0)=(Y, \Delta)$. 
\item[(2)] Let $U$ be the largest open subset of $Y$ such that 
$(U, D|_U)$ is a simple normal crossing pair. 
Then $\pi_i$ is an isomorphism over $U$ for every $i\geq 1$. 
\item[(3)] $D_i$ is the strict transform of $D$ on $Y_i$ for every $i\geq 1$. 
\item[(4)] $\pi_i$ is a blow-up whose center is a stratum of $(Y_{i-1}, 
\Supp \Delta_{i-1})$ and is located outside $U$ for every $i\geq 1$. 
\item[(5)] $\Delta_i=(\pi_i)^{-1}_*\Delta_{i-1}+E_i$, where 
$E_i$ is the $\pi_i$-exceptional divisor on $Y_i$ for every $i\geq 1$. 
\item[(6)] $(Y_k, D_k)$ is a simple normal crossing pair, that is, 
$D_k$ is Cartier {\em{(}}resp.~$\mathbb Q$-Cartier or 
$\mathbb R$-Cartier{\em{)}}. 
\end{itemize}
We note that $(Y_i, \Delta_i)$ is a simple normal crossing 
pair for every $i\geq 1$. 
\end{lem}
\begin{proof}
This is a direct consequence of 
\cite[8 The non-reduced case]{bierstone}. 
\end{proof}
\begin{lem}\label{ff-lem2.12}
Let $Y$ be an equidimensional variety which satisfies Serre's 
$S_2$ condition and is simple normal crossing 
in codimension one. 
Let $\Delta_Y$ be an effective $\mathbb Q$-divisor 
{\em{(}}resp.~$\mathbb R$-divisor{\em{)}} 
on $Y$ such that 
no irreducible component of 
$\Delta_Y$ is contained in the singular locus of $Y$. 
Assume that there exists a non-empty dense Zariski open set $Y_0$ 
of $Y$ such that $(Y_0, \Delta_{Y_0})$ 
is semi-log-canonical, where 
$\Delta_{Y_0}=(\Delta_Y)|_{Y_0}$. 
Then there exists a projective surjective 
birational morphism 
$\pi:Y'\to Y$, which is a composite of 
blow-ups, from a simple normal crossing variety $Y'$ with the following 
properties. 
\begin{itemize}
\item[(1)] $\Sing Y'$ maps birationally onto the 
closure of $\Sing Y^{\mathrm{snc}}$ by $\pi$, 
where $Y^{\mathrm{snc}}$ is the open subset of $Y$ which 
has only smooth points and simple normal crossing points. 
Note that $\Sing Y'$ {\em{(}}resp.~$\Sing Y^{\mathrm{snc}}${\em{)}} 
is the singular 
locus of $Y'$ {\em{(}}resp.~$Y^{\mathrm{snc}}${\em{)}}. 
\item[(2)] $\Delta_{Y'}$ is a subboundary $\mathbb Q$-divisor 
{\em{(}}resp.~$\mathbb R$-divisor{\em{)}} on 
$Y'$, that is, $(\Delta_{Y'})^{\leq 1}=\Delta_{Y'}$, 
such that $(Y', \Delta_{Y'})$ is a simple normal crossing 
pair with $K_{Y'_0}+\Delta_{Y'_0}=\pi^*_0(K_{Y_0}+\Delta_{Y_0})$, 
where $Y'_0=\pi^{-1}(Y_0)$, 
$\Delta_{Y'_0}=(\Delta_{Y'})|_{Y'_0}$, and $\pi_0=\pi|_{Y'_0}$. 
\item[(3)] $\Delta_{Y'}$ is the closure of $\Delta_{Y'_0}$ in $Y'$. 
\item[(4)] $\pi$ is an isomorphism over $U$, where 
$U$ is the largest open subset of $Y_0$ such that 
$(U, (\Delta_Y)|_U)$ is a simple normal crossing pair. 
\end{itemize}
\end{lem}
\begin{proof}
By \cite[Theorem 1.4]{bierstone}, 
we can construct a projective surjective birational morphism 
$\pi^\dag: Y^\dag\to Y$, which is a composite of 
blow-ups, 
from a simple normal crossing 
variety $Y^\dag$ such that 
$\Sing Y^\dag$ maps birationally onto the closure 
of $\Sing Y^{\mathrm{snc}}$ and 
that $(Y^\dag, (\pi^\dag)^{-1}_*\Delta_Y+E^\dag)$ is a simple 
normal crossing pair. We note that 
$E^\dag$ is a reduced Weil divisor on $Y^\dag$ whose support coincides with the 
exceptional locus of $\pi^\dag$. Of course, $\pi^\dag$ is 
an isomorphism over $\widetilde U(\supset U)$ by construction, 
where $\widetilde U$ is the largest open subset 
of $Y$ such that 
$(\widetilde U, (\Delta_Y)|_{\widetilde U})$ is a simple normal crossing 
pair. 
We put 
$$
K_{Y^\dag_0}+\Delta_{Y^\dag_0}=
(\pi^\dag_0)^*(K_{Y_0}+\Delta_{Y_0}), 
$$ 
where $Y^\dag_0=(\pi^\dag)^{-1}(Y_0)$ and $\pi^\dag_0=\pi^\dag|_{Y_0}$. 
Note that $\Delta_{Y^\dag_0}$ is a subboundary $\mathbb Q$-divisor 
(resp.~$\mathbb R$-divisor) on $Y^\dag_0$ since $(Y_0, \Delta_{Y_0})$ is 
semi-log-canonical. 
By construction, we see that 
$\Supp \Delta_{Y^\dag_0}\subset \Supp ((\pi^\dag)^{-1}_*
\Delta_Y+E^\dag)$. 
Then the closure of $\Delta_{Y^\dag_0}$ is contained 
in a simple normal crossing divisor $\Supp ((\pi^\dag)^{-1}_*
\Delta_Y+E^\dag)$. 
By Lemma \ref{ff-lem2.11}, 
we can take some sequence of blow-ups $\pi':Y'\to Y^\dag$, which is an 
isomorphism over $Y^\dag_0$ and is an isomorphism 
over the generic point of every stratum of $Y^\dag$, 
the closure $\Delta_{Y'}$ of $\Delta_{Y'_0}=\Delta_{Y^\dag_0}$ 
is a subboundary $\mathbb Q$-divisor (resp.~$\mathbb R$-divisor) 
on $Y'$ such that 
$(Y', \Delta_{Y'})$ is a simple normal crossing pair. 
By construction, we see that this is what we wanted. 
\end{proof}

Let us consider blow-ups of simple normal crossing varieties. 

\begin{lem}\label{ff-lem2.13}
Let $Y$ be a projective simple normal crossing variety. 
Let $\pi:Y'\to Y$ be a blow-up whose center $S$ is a stratum of 
$Y$ with $\codim _Y S\geq 2$. 
We take a Weil divisor $K_{Y'}$ on $Y'$ such that 
$\omega_{Y'}\simeq \mathcal O_{Y'}(K_{Y'})$ and that 
no irreducible component of $K_{Y'}$ is 
contained in the singular locus of $Y'$. 
Then $\omega_Y\simeq \mathcal O_Y(K_Y)$ and 
$K_{Y'}+E=\pi^*K_Y$ hold, where $K_Y=\pi_*K_{Y'}$ and $E=\Exc(\pi)$, 
that is, $E$ is a reduced Weil divisor on $Y'$ whose support 
coincides with the exceptional locus of $\pi$. 
Let $\Delta$ be a Cartier divisor on $Y$. 
Then the inclusion 
$$\pi_*\mathcal O_{Y'}(kK_{Y'}+\Delta')\subset\mathcal O_Y
(kK_Y+\Delta)$$ holds for every positive integer $k$, 
where $\Delta'=\pi^*\Delta$. 
\end{lem}
\begin{proof}
Since $\codim _Y S\geq 2$, the natural inclusion 
$\mathcal O_Y\hookrightarrow \pi_*\mathcal O_{Y'}$ 
is an isomorphism. 
We put $K_Y=\pi_*K_{Y'}$. 
Then $K_Y$ satisfies $\omega_Y\simeq \mathcal O_Y(K_Y)$. 
By definition, $K_{Y'}-\pi^*K_Y$ is contained in the exceptional 
locus of $\pi$. 
By considering the normalizations of $Y'$ and $Y$, 
we see that $K_{Y'}+E=\pi^*K_Y$. 
Since $\mathcal O_Y\simeq \pi_*\mathcal O_{Y'}$, we get the desired 
inclusion 
$$
\pi_*\mathcal O_{Y'}(kK_{Y'}+\Delta')\subset 
\pi_*\mathcal O_{Y'}(k(K_{Y'}+E)+\Delta') 
\simeq \mathcal O_Y(kK_Y+\Delta)
$$ 
for every positive integer $k$ since $E$ is an 
effective divisor on $Y'$. 
\end{proof}
\begin{lem}\label{ff-lem2.14}
Let $Y$ be a projective simple normal crossing variety. 
Let $\pi:Y'\to Y$ be a blow-up whose center $S$ is a minimal 
stratum of 
$Y$ with $\codim _Y S=1$. 
We take a Weil divisor $K_Y$ on $Y$ such that 
$\omega_Y\simeq \mathcal O_Y(K_Y)$ and that 
no irreducible component of $K_Y$ is contained in the 
singular locus of $Y$. Then $\omega_{Y'}\simeq 
\mathcal O_{Y'}(K_{Y'})$ holds, where 
$K_{Y'}=\pi^*K_Y-E$ and $E=\Exc (\pi)$, that is, $E$ is a 
reduced Weil divisor on $Y'$ whose support coincides with 
the exceptional locus of $\pi$. 
Let $\Delta$ be a Cartier divisor on $Y$. 
Then the inclusion 
$$\pi_*\mathcal O_{Y'}(kK_{Y'}+\Delta')\subset\mathcal O_Y
(kK_Y+\Delta)$$ holds for every positive integer $k$, 
where $\Delta'= \pi^*\Delta$. 
\end{lem}
\begin{proof}
By assumption, $S$ is a minimal stratum of $Y$ with $\codim _Y S=1$. 
This means that there are two irreducible components 
$Y_1$ and $Y_2$ of $Y$ such that $S$ is 
an irreducible component of $Y_1\cap Y_2$ and 
that $S\cap Y_i=\emptyset$ for every irreducible component $Y_i$ 
other than $Y_1$ and $Y_2$. 
Therefore, $\pi:Y'\to Y$ is nothing but the normalization in 
a neighborhood of $S$. 
Thus we see that $K_{Y'}=\pi^*K_Y-E$ satisfies $\omega_{Y'}
\simeq \mathcal O_{Y'}(K_{Y'})$. 
By considering the dual of the natural inclusion $\mathcal O_Y
\hookrightarrow \pi_*\mathcal O_{Y'}$, we get a generically 
isomorphic injection $\pi_*\omega_{Y'}\hookrightarrow 
\omega_Y$. 
By projection formula and $\pi^*K_Y=K_{Y'}+E$, we obtain 
the desired inclusion 
\begin{equation*}
\begin{split}
\pi_*\mathcal O_{Y'}(kK_{Y'}+\Delta')&\subset \pi_*\mathcal O_{Y'} 
(K_{Y'}+(k-1)\pi^*K_Y+\pi^*\Delta) \\ 
&\simeq \pi_*\mathcal O_{Y'}(K_{Y'})\otimes 
\mathcal O_Y((k-1)K_Y+\Delta)\\ 
&\subset \mathcal O_Y(kK_Y+\Delta)
\end{split} 
\end{equation*}
for every positive integer $k$ since $\pi_*\mathcal O_{Y'}(K_{Y'})\subset 
\mathcal O_Y(K_Y)$ as mentioned above. 
\end{proof}
\end{say}

\section{A quick review of $\mathcal M^{stable}$}\label{ff-sec3} 

In this section, we quickly review the moduli space of stable varieties 
(see \cite{kollar2}). 
First, let us recall the definition of {\em{stable varieties}}. 

\begin{defn}[Stable varieties]\label{ff-def3.1} 
Let $X$ be a connected projective semi-log-canonical 
variety with ample canonical divisor. 
Then $X$ is called a {\em{stable variety}} or a {\em{semi-log-canonical model}}. 
\end{defn}

In order to obtain the boundedness of the moduli functor of stable varieties, we 
have to 
fix some numerical invariants. So we introduce the 
notion of the {\em{Hilbert function}} for stable 
varieties. 

\begin{defn}[Hilbert function of stable varieties]\label{ff-def3.2} 
Let $X$ be a stable variety. The {\em{Hilbert function}} of 
$X$ is 
$$
H_X(m):=\chi (X, \omega_X^{[m]})
$$ 
where $\omega_X^{[m]}=
(\omega_X^{\otimes m})^{**}\simeq \mathcal O_X(mK_X)$. 
By \cite[Corollary 1.9]{fujino}, we see that 
$$
H_X(m)=\dim _{\mathbb C} H^0(X, \mathcal O_X(mK_X))\geq 0
$$ 
for every $m\geq 2$. 
\end{defn}

The following definition of the 
moduli functor of stable varieties is mainly due to Koll\'ar. 
Note that a stable variety $X$ is not necessarily Cohen--Macaulay when 
$\dim X\geq 3$. 
We think that it is one of the main difficulties when we treat families of 
stable varieties. 

\begin{defn}
[Moduli functor of stable varieties]\label{ff-def3.3} 
Let $H(m)$ be a $\mathbb Z$-valued function. 
The moduli functor of stable varieties with Hilbert function $H$ is 
$$
\mathcal M^{stable}_H(S):=
\left\{
\begin{array}{c}
\mbox{Flat, proper families $X\to S$, fibers are stable}\\ 
\mbox{varieties with ample canonical divisor}\\
\mbox{and Hilbert function $H(m)$, $\omega^{[m]}_{X/S}$ is flat over $S$}\\ 
\mbox{and commutes with base change for every $m$, }\\
\mbox{modulo isomorphisms over $S$.}
\end{array}
\right\}.
$$
\end{defn} 

\begin{rem}\label{ff-rem3.4}
We consider $(f:X\to S)\in \mathcal M^{stable}_H(S)$. 
By the base change theorem and \cite[Corollary 1.9]{fujino}, 
we obtain that $f_*\omega^{[m]}_{X/S}$ is a locally free sheaf on 
$S$ with $\mathrm{rank} f_*\omega^{[m]}_{X/S}=H(m)$ for every $m\geq 2$. 
\end{rem}

Let us quickly review the construction of the coarse moduli space of 
stable varieties following \cite{kollar2}. 

\begin{say}[Coarse moduli space of $\mathcal M^{stable}$]\label{ff-say3.5}
Let us consider the moduli functor  
$$
\mathcal M^{stable}(S):=
\left\{
\begin{array}{c}
\mbox{Flat, proper families $X\to S$, fibers are stable}\\ 
\mbox{varieties with ample canonical divisor, }\\
\mbox{$\omega^{[m]}_{X/S}$ is flat over $S$}\\ 
\mbox{and commutes with base change for every $m$, }\\
\mbox{modulo isomorphisms over $S$.}
\end{array}
\right\}
$$ 
of stable varieties. It is obvious that $\mathcal M^{stable}_H$ is an open and 
closed subfunctor of $\mathcal M^{stable}$. 
It is known that the moduli functor $\mathcal M^{stable}$ is well-behaved, that is, 
$\mathcal M^{stable}$ is locally closed.  
For the details, see \cite[Corollary 25]{kollar-hull}. 
We have already known that the moduli functor $\mathcal M^{stable}$ 
satisfies the valuative criterion of separatedness and the valuative criterion of 
properness by Koll\'ar's gluing theory and the existence of log canonical closures 
(see \cite[Theorem 24]{kollar-s} and \cite[Section 7]{hacon-xu}). 
Moreover, it is well known that the automorphism group $\Aut (X)$ of 
a stable variety $X$ is a finite group (for a more general 
result, see \cite[Corollary 6.17]{fujino}). 
Then, by using \cite[1.2 Corollary]{keel-mori}, 
we obtain a coarse moduli space $M^{stable}$ of 
$\mathcal M^{stable}$ in the category of algebraic spaces 
(see, for example, \cite[5.3 (Existence of coarse moduli spaces)]{kollar2}). 
Note that 
$M^{stable}$ is a separated algebraic space which is locally of finite type. 
Since $\mathcal M^{stable}$ satisfies the valuative 
criterion of properness, $M^{stable}_H$ 
is proper if and only if it is of finite type. 
\end{say}

\section{Proof of theorems}\label{ff-sec4}

Let us start the proof of Theorem \ref{ff-thm1.9}.
Theorem \ref{ff-thm1.9} is essentially contained in \cite[Section 7]{fujino-fujisawa} 
(see also \cite{ffs}).
We need no extra assumptions on $D$ and local monodromies since $C$ is a curve.

\begin{proof}[Proof of Theorem \ref{ff-thm1.9}]
In Step \ref{ff-step1.10-1}, we will reduce the 
problem to a simpler case by using 
\cite{bierstone}. In Step \ref{ff-step1.10-2}, we will prove the desired 
nefness by using \cite{fujino-fujisawa}. 
\begin{step}\label{ff-step1.10-1}
There is a closed subset $\Sigma$ of $C$ such that
every stratum of $(X, D)$ is smooth over $C_0=C\setminus \Sigma$. 
This means that every stratum of $(X_0, D_0)$ is smooth 
over $C_0$, 
where $X_0=f^{-1}(C_0)$ and 
$D_0=D|_{X_0}$. 
Apply \cite[Theorem 1.4]{bierstone} to
$(X, \Supp (D+f^*\Sigma))$. Then we
obtain a birational morphism $g:X'\to X$ from a projective simple normal crossing
variety $X'$ such that $g$ is an isomorphism
outside $\Supp f^*\Sigma$ and
that $g_*^{-1}D+g^*f^*\Sigma$ has a simple normal crossing support on $X'$.
Let $D'$ be the horizontal part of
$g_*^{-1}D$. 
By taking some more blow-ups, we may further 
assume that $D'$ is a Cartier divisor on $X'$ (see 
Lemma \ref{ff-lem2.11}). We note that 
$g:X'\to X$ is an isomorphism over $X_0$ by construction. 
Therefore, $(X', D')$ is a simple normal crossing
pair. 
By construction, we have $\omega_{X'}(D')\simeq 
g^*\omega_X(D)\otimes \mathcal O_{X'}(E)$ such that 
$f\circ g(E)\subset \Sigma$ and that 
the effective part of $E$ is $g$-exceptional. 
Thus, we have $g_*\omega_{X'}(D')\subset \omega_X(D)$. 
Hence we obtain an inclusion 
$$
f_*g_*\omega_{X'/C}(D')\to f_*\omega_{X/C}(D),
$$
which is an isomorphism over $C_0$. 
Therefore, it is sufficient to
prove that $f_*g_*\omega_{X'/C}(D')$ is nef by Lemma \ref{ff-lem2.2}.
By replacing $(X, D)$ with $(X', D')$, we may assume that
every stratum of $(X, D)$ is dominant onto $C$. 
Of course, by assumption, every stratum of $(X_0, D_0)$ is smooth over 
$C_0$. 
\end{step}
\begin{step}\label{ff-step1.10-2}
By \cite[Theorem 4.15]{fujino-fujisawa}, 
$R^df_{0*}\iota_!\mathbb Q_{X_0\setminus D_0}$, where 
$d=\dim X-1$, $X_0=f^{-1}(C_0)$,
$f_0=f|_{X_0}$,
$D_0=D|_{X_0}$, and
$\iota:X_0\setminus D_0\hookrightarrow X_0$, underlines 
a graded polarizable admissible variation of $\mathbb Q$-mixed 
Hodge structure. 
In particular, every local monodromy on $R^df_{0*}\iota_!
\mathbb Q_{X_0\setminus D_0}$ around $\Sigma$ is quasi-unipotent 
(see \cite[Definition 3.11]{fujino-fujisawa}). 
Moreover, we can consider (upper and lower) 
canonical extensions of Hodge bundles (see \cite[Definition 3.11 and 
Remark 7.4]{fujino-fujisawa}). By \cite[Theorem 7.3 (a)]{fujino-fujisawa}, 
$R^df_*\mathcal O_X(-D)$ is characterized as the lower 
canonical extension of 
$$
\Gr^0_F(R^df_{0*}\iota_!\mathbb Q_{X_0\setminus D_0}\otimes 
\mathcal O_{C_0}). 
$$ 
We note that we can freely replace $C_0$ with its non-empty Zariski 
open set. We take a unipotent reduction $\pi:C'\to C$ 
of $R^df_{0*}\iota _!\mathbb Q_{X_0\setminus D_0}$ (see \cite[Theorem 17 and 
Corollary 18]{kawamata-abel}).
We may assume that $\pi$ is a finite Galois cover 
(see \cite[Theorem 1-1-1]{kmm}). By shrinking $C_0$, we may further
assume that
$\pi:C'\to C$ is \'etale over $C_0$. 
We note that the local system 
$\pi_0^*R^df_{0*}\iota_!\mathbb Q_{X_0\setminus D_0}$ 
on $C'_0=\pi^{-1}(C_0)$ underlies a 
graded polarizable admissible variation of $\mathbb Q$-mixed Hodge structure 
since $R^df_{0*}\iota_!\mathbb Q_{X_0\setminus D_0}$ underlies 
a graded polarizable admissible variation of 
$\mathbb Q$-mixed Hodge structure, 
where
$\pi_0=\pi|_{C_0'}:C_0'\to C_0$.
Let $\mathcal G$ be the canonical
extension of
$$\Gr_F^0(\pi_0^*R^df_{0*}\iota_!\mathbb Q_{X_0\setminus D_0}
\otimes \mathcal O_{C_0'}). $$  
Then $\mathcal G$ is locally free and $\mathcal G^*$ is a 
nef locally free sheaf on $C'$ 
(see \cite[Corollary 5.23]{fujino-fujisawa}, 
\cite{ffs}, \cite{fujino-fujisawa2}, 
\cite{fujisawa}, and so on). 
Since $R^df_*\mathcal O_X(-D)\simeq (\pi_*\mathcal G)^G$,
where $G$ is the Galois group of $\pi:C'\to C$, we obtain a nontrivial
map
$\pi^*R^df_*\mathcal O_X(-D)\to \mathcal G$,
which is an isomorphism on $C_0'$.
Note that $R^df_*\mathcal O_X(-D)$ is the lower canonical
extension of
$$
\Gr_F^0(R^df_{0*}\iota_!\mathbb Q_{X_0\setminus D_0}
\otimes \mathcal O_{C_0}).
$$
Therefore, by taking the dual,
we obtain an inclusion $0\to \mathcal G^*\to \pi^*f_*\omega_{X/C}(D)$,
which is an isomorphism on $C_0'$. Thus,
$\pi^*f_*\omega_{X/C}(D)$ is nef by Lemma \ref{ff-lem2.2}.
So we obtain that $f_*\omega_{X/C}(D)$ is nef because $\pi$ is
surjective.
\end{step}
Anyway, we obtain the desired nefness of $f_*\omega_{X/C}(D)$. 
\end{proof}

For an alternative proof of Theorem \ref{ff-thm1.9} 
based on the Koll\'ar--Ohsawa type vanishing theorem for semi-log-canonical 
pairs, see \cite{fujino-new}. 
We note that \cite{fujino-new} depends on the theory of 
mixed Hodge structures on cohomology with compact support. 

\begin{rem}\label{ff-rem4.1}
When $X$ is smooth in Theorem \ref{ff-thm1.9},
the semipositivity theorem
obtained in \cite[Theorem 3.9]{fujino-high}
is sufficient for the proof of Theorem \ref{ff-thm1.9}.
Note that \cite[Theorem 3.9]{fujino-high} also follows from the
theory of graded polarizable admissible variations of mixed Hodge structure.
\end{rem}

Let us prove Theorem \ref{ff-thm1.10} since it contains Theorem \ref{ff-thm1.7} 
as a special case. 

\begin{proof}[Proof of Theorem \ref{ff-thm1.10}]
In Step \ref{ff-step1.11-1}, we will prove the nefness of $f_*\omega_{X/C}(D)$. 
In Step \ref{ff-step1.11-2}, we will treat $f_*\mathcal O_X(m(K_{X/C}+D))$. 
\setcounter{step}{0}
\begin{step}\label{ff-step1.11-1}
We take a double cover $\pi:\widetilde X\to X$ due to Koll\'ar (see \cite[5.23]{kollar-book}). 
Then $\omega_X(D)$ is a direct summand of $\pi_*\omega_{\widetilde X}(\pi^{-1}D)$. 
By replacing $X$ and $f$ with $\widetilde X$ and $f\circ \pi$, respectively, 
we may further assume that $X$ is simple normal crossing in codimension one. 

We note that $X$ and $D$ satisfy the assumptions in Lemma \ref{ff-lem2.12} and 
that $(f^{-1}(U), D|_{f^{-1}(U)})$ is semi-log-canonical by assumption. 
Therefore, we can apply Lemma \ref{ff-lem2.12}. 
Then we can construct a projective surjective birational morphism 
$h:Z\to X$, which is a composite of blow-ups, from 
a simple normal crossing variety $Z$ with the following properties. 
\begin{itemize}
\item[(1)] there exists a subboundary $\mathbb Q$-divisor 
$B$ on $Z$, that is, $B^{\leq 1}=B$, 
such that $(Z, B)$ is a simple normal crossing pair. 
\item[(2)] $\Sing Z$ maps birationally onto 
the closure of $\Sing X^{\mathrm{snc}}$ by $h$, where 
$\Sing Z$ (resp.~$\Sing X^{\mathrm{snc}}$) is the singular locus of 
$Z$ (resp.~$X^{\mathrm{snc}}$). 
Note that $X^{\mathrm{snc}}$ is the open subset of $X$ 
which has only smooth points and simple normal crossing 
points. 
\item[(3)] $K_V+B|_V=g^*\left(K_{f^{-1}(U)}+D|_{f^{-1}(U)}\right)$, where 
$V=(f\circ h)^{-1}(U)$ and $g=h|_V: V\to f^{-1}(U)$. 
\end{itemize}
We note that the natural inclusion $\mathcal O_X\hookrightarrow 
h_*\mathcal O_Z$ is an isomorphism and that 
$Z\setminus V$ contains no irreducible components of $B$ by construction. 
Thus we have 
$$
g_*\omega_V((B|_V)^{=1})\simeq 
\omega_{f^{-1}(U)}(D|_{f^{-1}(U)}) \quad  \text{and}\quad  
h_*\omega_Z(B^{=1})\subset \omega_X(D). 
$$ 
Therefore, it is sufficient to prove that 
$(f\circ h)_*\omega_{Z/C}(B^{=1})$ is nef by Lemma \ref{ff-lem2.2} 
in order to prove the nefness of $f_*\omega_{X/C}(D)$. 
By shrinking $U$ suitably, we may assume that every stratum of $Z$ maps 
to a point in $C\setminus U$ or is dominant onto $C$. 
If $Z\setminus V$ contains an irreducible component of $B$, then we take some 
blow-ups outside $V$ and replace $B$ with 
the closure of $B|_V$ (see Lemma \ref{ff-lem2.11}). 
Thus, we can always 
assume that $Z\setminus V$ contains no irreducible components of $B$. 

Assume that there exists 
a stratum $S$ of $Z$ in $(f\circ h)^{-1}(\Sigma)$, where 
$\Sigma =C\setminus U$. 
If $\codim _Z S\geq 2$ or $S$ is a minimal stratum of 
$Z$ with $\codim _Z S=1$, then we take a blow-up $\pi:Z'\to Z$ of 
$Z$ along $S$ as in Lemmas \ref{ff-lem2.13} and \ref{ff-lem2.14}. 
We note that $\pi^*B=\pi^{-1}_*B$ and that 
$(\pi^*B)^{=1}=\pi^{-1}_*(B^{=1})=\pi^*(B^{=1})$. 
By Lemmas \ref{ff-lem2.13} and \ref{ff-lem2.14}, there is a 
generically isomorphic inclusion 
$$
\pi_*\omega_{Z'}((\pi^*B)^{=1})\subset 
\omega_Z(B^{=1}), 
$$ 
which is an isomorphism over $U$. 
Therefore, by Lemma \ref{ff-lem2.2}, it is sufficient to prove that 
$(f\circ h\circ \pi)_*\omega_{Z'/C}((\pi^*B)^{=1})$ is nef. 
This means that we can replace $(Z, B)$ with $(Z', \pi^*B)$. 
By repeating this process finitely many times, we may assume that 
$(f\circ h)^{-1}(\Sigma)$ contains no strata of $Z$. In this case, 
the nefness of $(f\circ h)_*\omega_{Z/C}(B^{=1})$ follows from 
Theorem \ref{ff-thm1.9}. Anyway, we obtain that $f_*\omega_{X/C}(D)$ is 
a nef locally free sheaf on $C$. 
\end{step}

\begin{step}[{see Proof of \cite[Corollary 2.45]{viehweg2}}]\label{ff-step1.11-2}
By Viehweg's clever covering trick,
we can prove that $f_*\mathcal O_X(m(K_{X/C}+D))$ is nef 
for every $m\geq 2$ by using the case where $m=1$. 
Here, we closely follow the proof of \cite[Corollary 2.45]{viehweg2}.

Let $\mathcal H$ be an ample line bundle on $C$.
We set
$$
r=\min \left\{\mu\in \mathbb Z_{>0}\,\left|
\,(f_*\mathcal O_X(k(K_{X/C}+D)))\otimes 
\mathcal H^{\mu k-1}\, \text{is nef}\right.\right\}.
$$
By 
assumption, 
we have that the natural map 
$$f^*f_*\mathcal O_X(k(K_{X/C}+D))\to \mathcal O_X(k(K_{X/C}+D))$$ 
is surjective. 
Since 
$(f_*\mathcal O_X(k(K_{X/C}+D)))\otimes \mathcal H^{rk-1}$ is a nef locally free sheaf, 
$(f_*\mathcal O_X(k(K_{X/C}+D)))\otimes \mathcal H^{rk}$ is ample. 
Therefore, we see that 
$S^N((f_*\mathcal O_X(k(K_{X/C}+D)))\otimes \mathcal H^{rk})$ is generated by its 
global sections 
for some positive integer $N$. 
Hence we see that $\mathcal O_X(K_{X/C}+D)\otimes f^*\mathcal  H^r$ is semi-ample.
More precisely,
$\mathcal O_X(k(K_{X/C}+D))\otimes f^*\mathcal H^{rk}$ is locally free and semi-ample. 
By the usual covering
argument (see Remark \ref{ff-rem4.2}),
$(f_*\mathcal O_X(k(K_{X/C}+D)))\otimes \mathcal H^{r(k-1)}$ is nef (see 
\cite[Proposition 2.43]{viehweg2}).
This is only possible if $(r-1)k-1<r(k-1)$.
It is equivalent to $r\leq k$.
Therefore,
$(f_*\mathcal O_X(k(K_{X/C}+D)))\otimes \mathcal H^{k^2-1}$ is nef.
The same holds true if we take any base change by $\pi:C'\to C$ such that
$\pi$ is a finite morphism from a smooth 
projective curve and ramifies only over general points of $C$.
Therefore, $f_*\mathcal O_X(k(K_{X/C}+D))$ is nef 
(see Lemma \ref{ff-lem2.3}).
By the same argument as above, we 
see that $\mathcal O_X(K_{X/C}+D)\otimes f^*\mathcal H$ is semi-ample 
since $f_*\mathcal O_X(k(K_{X/C}+D))$ is nef. 
More precisely, $\mathcal O_X(k(K_{X/C}+D))\otimes f^*\mathcal H^k$ 
is locally free and semi-ample in the usual 
sense.
By the covering argument (see Remark \ref{ff-rem4.2}),
$(f_*\mathcal O_X(m(K_{X/C}+D)))\otimes \mathcal H^{m-1}$ is nef for
every $m>0$ (see \cite[Proposition 2.43]{viehweg2}).
The same holds true if we take any base change by $\pi:C'\to C$
such that $\pi$ is a finite 
morphism from a smooth projective curve and ramifies only over 
general points of $C$. 
Therefore,
$f_*\mathcal O_X(m(K_{X/C}+D))$ is nef for
every $m>0$ (see Lemma \ref{ff-lem2.3}). 
\end{step}
We have completed the proof of Theorem \ref{ff-thm1.10}.
\end{proof}

\begin{rem}[{see \cite[4.15.~Lemma and 4.16]{kollar}}]\label{ff-rem4.2}
Let $\varphi:X'\to X$ be a cyclic cover associated to a general
member $A\in |\mathcal O_X(kl(K_{X/C}+D))\otimes f^*
\mathcal H^{rkl}|$
(resp.~$|\mathcal O_X(kl(K_{X/C}+D))\otimes f^*\mathcal H^{kl}|$)
for some positive integer $l$.
Then $f':=f\circ \varphi:X'\to C$ and 
$\varphi^{-1}D$ 
satisfy all the assumptions for $f:X\to C$ and $D$.
Therefore,
we see that $f'_*\omega_{X'/C}(\varphi^{-1}D)$ is nef by Step \ref{ff-step1.11-1} in the 
proof of Theorem \ref{ff-thm1.10}.
By construction, 
 $\mathcal O_X(k(K_{X/C}+D))\otimes f^*\mathcal H^{r(k-1)}$
 (resp.~$\mathcal O_X(m(K_{X/C}+D))\otimes f^*\mathcal H^{m-1}$)
 is a direct summand of $\varphi_*\omega_{X'/C}(\varphi^{-1}D)$. 
 Thus, we obtain that $f_*\mathcal O_X(k(K_{X/C}+D))\otimes 
 \mathcal H^{r(k-1)}$ (resp.~$f_*\mathcal 
 O_X(m(K_{X/C}+D))\otimes \mathcal H^{m-1}$) is a nef locally free 
 sheaf on $C$. 
\end{rem}

Theorem \ref{ff-thm1.2} is almost obvious by Theorem \ref{ff-thm1.7}. 

\begin{proof}[Proof of Theorem \ref{ff-thm1.2}] 
We consider $(f:X\to C)\in \mathcal M^{stable}(C)$ where 
$C$ is a smooth projective curve. 
By Kawakita's inversion of adjunction \cite[Theorem]{kawakita} (see also 
\cite[Lemma 2.10 and Corollary 2.11]{patakfalvi}), 
$X$ itself is a semi-log-canonical variety. 
By the definition of $\mathcal M^{stable}$, 
we can find a positive integer $k$ such that $\omega^{[k]}_{X/C}$ is locally free and 
$f$-ample. Hence $f_*\omega^{[m]}_{X/C}$ is nef 
for every $m\geq 1$ by Theorem 
\ref{ff-thm1.7}. It implies that $\mathcal M^{stable}$ is semipositive in the sense of 
Koll\'ar. 
\end{proof}

By Koll\'ar's results (see \cite[Sections 2 and 3]{kollar2}), 
Theorem \ref{ff-thm1.1} follows from Theorem \ref{ff-thm1.2}. 
For some technical details, 
see also \cite[Theorem 4.34 and Theorem 9.25]{viehweg2}.  

\begin{proof}[Proof of Theorem \ref{ff-thm1.1}] 
It is sufficient to prove this theorem for connected subspaces. 
Let $Z$ be a connected closed complete subspace of $M^{stable}$. 
It is obvious that $M^{stable}$ has an open subspace of finite type which contains $Z$. 
By replacing $\mathcal M^{stable}$ with the subfunctor given by 
this subspace, we get a new functor $\mathcal N$ which is bounded. 
By recalling 
the construction of the coarse moduli space, 
we know that there is a locally closed subscheme $S$ of $\mathrm {Hilb}(\mathbb P^N)$ for some 
$N$ such that $Z$ is obtained as the geometric quotient $S/\Aut (\mathbb P^N)$. 
Let $f:X\to S$ be the universal family. By the proof of 
\cite[2.6.~Theorem]{kollar}, we see that 
$\det (f_*\omega^{[k]}_{X/S})^p$ descends to an ample line bundle 
on $Z$ for a sufficiently large and divisible 
positive integer $k$ and a sufficiently divisible positive integer $p$ 
(see \cite[Lemma 9.26]{viehweg2}). 
Note that \cite[2.6.~Theorem]{kollar} needs the semipositivity of $\mathcal M^{stable}$. 
For the details, see \cite[Sections 2 and 3]{kollar} and \cite{viehweg2}. 
\end{proof}

Theorem \ref{ff-thm1.10} is useful for the projectivity of the 
moduli space of {\em{stable maps}} 
(see \cite{fp}). For some related topics, see also \cite{alexeev2} and 
\cite{kp}. 

\begin{say}[{Projectivity of the space of stable maps (see \cite{fp})}]\label{ff-say4.3} 
We freely use the notation in \cite{fp}. 
Let $\mathcal F=(\pi, \mathcal C\to S, \{p_i\}, \mu)$ be a 
{\em{stable family of maps}} over $S$ to $\mathbb P^r$. 
For the definition, see \cite[1.1.~Definitions]{fp}. We set 
$$
E_k(\pi)=\pi_*\left(\omega_{\pi}^k\left(\sum _{i=1}^n kp_i\right)
\otimes \mu^*(\mathcal O(3k))
\right). 
$$ 
In \cite[Lemma 3]{fp}, it is proved that $E_k(\pi)$ is a nef 
locally free sheaf on $S$ for $k\geq 2$ 
by using the results in \cite[Section 4]{kollar}. 
This nefness 
is used for the projectivity of the moduli space of stable maps in \cite[4.3.~Projectivity]{fp}. 
The nefness of $E_k(\pi)$ can be checked 
as follows:  

Since $k\geq 2$, by the base change theorem, we may assume that 
$S$ is a smooth projective curve. 
We take a general member $H$ of $|\mu^*\mathcal O(3)|$. 
Then $(\mathcal C, \sum _{i=1}^np_i+H)$ is a semi-log-canonical 
surface and 
$K_{\mathcal C/S}+\sum _{i=1}^np_i+H$ is $\pi$-ample. Therefore 
$$\pi_*\mathcal O_{\mathcal C}\left(k\left(K_{\mathcal C/S}+\sum _{i=1}^np_i+H
\right)\right)\simeq 
E_k(\pi)$$ is nef for every $k\geq 2$ by Theorem \ref{ff-thm1.10}. 
\end{say}

Let us start the proof of Theorem \ref{ff-thm1.11}. 

\begin{proof}[Proof of Theorem \ref{ff-thm1.11}] 
We will use Viehweg's covering arguments (see \cite{viehweg1}) 
and a special case of Theorem \ref{ff-thm1.10}. 

We take a double cover $\pi:\widetilde X\to X$ 
due to Koll\'ar (see \cite[5.23]{kollar-book}). 
We put $\widetilde \Delta=\pi^{-1}\Delta$. 
Then $\mathcal O_{\widetilde X} (k(K_{\widetilde X}+\widetilde \Delta))$ is 
locally free and $f\circ \pi$-generated. 
Moreover, 
$\mathcal O_X(kl(K_{X/C}+\Delta))$ is a direct summand of 
$\pi_*\mathcal O_{\widetilde X}(kl(K_{\widetilde X}+\widetilde \Delta))$ for 
every $l\geq 1$. 
Therefore, by replacing $X$ and $\Delta$ with 
$\widetilde X$ and $\widetilde \Delta$, respectively, 
we may further assume that $\widetilde X$ is simple normal 
crossing in codimension one. 

Since $X$ and $\Delta$ satisfy the assumptions in Lemma \ref{ff-lem2.12}, 
we can apply Lemma \ref{ff-lem2.12}. 
Then there is a projective surjective birational morphism 
$h:Z\to X$, which is a composite of 
blow-ups, 
from a simple normal crossing variety $Z$ with the following 
properties. 
\begin{itemize}
\item[(1)] there is a subboundary $\mathbb Q$-divisor $\Delta_Z$ on $Z$, 
that is, $\Delta^{\leq 1}_Z=\Delta_Z$, such that 
$(Z, \Delta_Z)$ is a simple normal crossing pair. 
\item[(2)] $\Sing Z$ maps birationally 
onto the closure of $\Sing X^{\mathrm{snc}}$ 
by $h$, where $\Sing Z$ (resp.~$\Sing X^{\mathrm{snc}}$) is the 
singular locus of $Z$ (resp.~$X^{\mathrm{snc}}$). 
Note that $X^{\mathrm{snc}}$ is the open subset of $X$ which has 
only smooth points and simple normal crossing points. 
\item[(3)] $K_V+\Delta_V=(h_V)^*(K_{f^{-1}(U)}+\Delta|_{f^{-1}(U)})$, where 
$V=(f\circ h)^{-1}(U)$, $h_V=h|_V: V\to f^{-1}(U)$, and 
$\Delta_V=(\Delta_Z)|_V$. 
\end{itemize} 
Note that the natural inclusion $\mathcal O_X\hookrightarrow 
h_*\mathcal O_Z$ is an isomorphism and 
that $Z\setminus V$ contains no irreducible components of $\Delta_Z$ 
by construction. 
Therefore, we have $$(h_V)_*\mathcal O_V(k(K_V+\Delta_V^{> 0}))
\simeq \mathcal O_{f^{-1}(U)}(k(K_X+\Delta))$$ 
and 
$$h_*\mathcal O_Z(k(K_Z+
\Delta_Z^{>0}))\subset \mathcal O_X(k(K_X+\Delta)). $$ 
Hence it is sufficient to prove that $(f\circ h)_*\mathcal O_Z
(k(K_{Z/C}+\Delta_Z^{> 0}))$
is nef by Lemma \ref{ff-lem2.2} for the proof of the nefness of 
$f_*\mathcal O_X(k(K_{X/C}+\Delta))$. 
By shrinking $U$, we may assume that every stratum of $(Z, \Supp \Delta_Z)$ is 
smooth over $U$ or maps to a point in $C\setminus U$. 
If $Z\setminus V$ contains an irreducible component of 
$\Delta_Z$, then we take some blow-ups outside $V$ 
and replace $\Delta_Z$ with 
the closure of $\Delta_V$ (see Lemma \ref{ff-lem2.11}). 
Therefore, we can assume that $Z\setminus V$ contains no irreducible components of 
$\Delta_Z$. 
We may further assume that 
$Z\setminus V$ contains no strata of $(Z, \Supp \Delta_Z)$ contained 
in $\Supp \Delta_Z$ 
by taking some more blow-ups.

Assume that there exists a stratum $S$ of $Z$ in $(f\circ h)^{-1}(\Sigma)$, where 
$\Sigma =C\setminus U$. 
If $\codim _Z S\geq 2$ or $S$ is a minimal stratum of $Z$ with 
$\codim _Z S=1$, then we take a blow-up $\pi:Z'\to Z$ of $Z$ along $S$ as in Lemmas 
\ref{ff-lem2.13} and \ref{ff-lem2.14}. 
We note that $k\Delta^{>0}_Z$ is Cartier and 
$\pi^*\Delta_Z=\pi^{-1}_*\Delta_Z$. 
Of course, $\pi^*(k\Delta^{>0}_Z)=k(\pi^{-1}_*\Delta_Z)^{>0}$ is 
Cartier. 
By Lemmas \ref{ff-lem2.13} and \ref{ff-lem2.14}, 
we have a natural inclusion 
$$
\pi_*\mathcal O_{Z'}(k(K_{Z'}+(\pi^*\Delta_Z)^{>0}))\subset 
\mathcal O_Z(k(K_Z+\Delta^{>0}_Z)), 
$$ 
which is an isomorphism over $U$. 
Therefore, we can replace $(Z, \Delta_Z)$ with $(Z', \pi^*\Delta_Z)$ by Lemma 
\ref{ff-lem2.2} in order to prove the nefness of $(f\circ h)_*\mathcal O_Z(k(K_{Z/C}+\Delta_Z))$. 
By repeating this process finitely many times, we may assume that 
$Z\setminus V$ contains no strata of $Z$. 

\begin{lem}\label{ff-lem4.4} 
Let $\pi:Z'\to Z$ be a projective surjective birational morphism 
from a simple normal crossing variety $Z'$, which is a 
composite of blow-ups whose centers are outside $V$. 
Assume that $(Z', \pi^{-1}_*(\Delta^{>0}_Z))$ is a simple normal 
crossing pair. Then we can write 
$$
K_{Z'}+\pi^{-1}_*(\Delta^{>0}_Z)=\pi^*(K_Z+\Delta^{>0}_Z)+E'
$$ 
for some effective $\pi$-exceptional $\mathbb Q$-divisor 
$E'$ on $Z'$.  
\end{lem}
\begin{proof}[Proof of Lemma \ref{ff-lem4.4}]
This follows from the fact that $Z\setminus V$ contains 
no strata of $Z$ and no strata of $(Z, \Supp \Delta_Z)$. 
\end{proof}
By the above construction, we have 
$$
K_V+\Delta_{V}^{>0}=h^*(K_{f^{-1}(U)}+\Delta_{f^{-1}(U)})+(-\Delta_V^{<0}). 
$$ 
We will apply the covering argument discussed in \cite[Section 4.4]{campana} 
(see also \cite[Section 8]{fujino-weak}), which is a modification of 
Viehweg's covering argument in \cite[Lemma 5.1 and Corollary 5.2]{viehweg1}. 
We set $g=f\circ h$. 
By taking more blow-ups over 
$Z\setminus V$ (see \cite{bierstone}, Lemma \ref{ff-lem2.11}, and 
Lemma \ref{ff-lem4.4}), we may assume that 
$$
\mathcal F:=\mathrm{Image}\left(g^*g_*
\mathcal O_Z(k(K_{Z/C}+\Delta
_Z^{>0}))\to 
\mathcal O_Z(k(K_{Z/C}+\Delta_Z^{>0}))\right)
$$
is a line bundle which is $g$-generated and that 
$$\mathcal O_Z(k(K_{Z/C}+\Delta_Z^{>0}))\simeq 
\mathcal F\otimes \mathcal O_Z(E)$$ such that 
$E$ is an effective Cartier divisor 
on $Z$ and that $\Supp E$ is a simple normal 
crossing divisor 
on $Z$. We note that 
$E$ is equal to $-k\Delta^{<0}_V$ over $U$ by construction. 
We may further assume that $\Supp E\cup \Supp \Delta^{>0}_Z$ is a simple 
normal crossing divisor on $Z$. 
Let $$g: Z\overset{\psi}{\longrightarrow} \widetilde C\longrightarrow C$$ 
be the Stein factorization. 
Note that $g_*\mathcal O_Z$ is a torsion-free sheaf on $C$ since every irreducible 
component of $Z$ is dominant onto $C$. 
We also note that $\widetilde C$ is normal (see \cite[Lemma 3.6]{fujino}). 
Then 
\begin{equation}\label{ff-eq4.1}
\mathcal F=\mathrm{Image}\left(\psi^*\psi_*
\mathcal O_Z(k(K_{Z/C}+\Delta
_Z^{>0}))\to 
\mathcal O_Z(k(K_{Z/C}+\Delta_Z^{>0}))\right). 
\end{equation}  
Let $\mathcal H$ be an ample line bundle on $C$. 
We set 
$$
r=\min \left\{\mu \in \mathbb Z_{>0}\, 
\left|\, g_*\mathcal O_Z(k(K_{Z/C}+\Delta_{Z}^{>0}))
\otimes \mathcal H^{\mu k-1} \ \text{is nef}\right.\right\}. 
$$

The following lemma is essentially 
contained in \cite[Lemma 5.1 and Corollary 5.2]{viehweg1}. 
 
\begin{lem}[{see \cite[Lemma 4.19]{campana} and \cite[Lemma 8.2]{fujino-weak}}]
\label{ff-lem4.5}
Let $g: Z\to C$ be as above. 
Let $\mathcal A$ be an ample line bundle on $C$. 
Assume that 
$$
S^N\left(g_*\mathcal 
O_Z(k(K_{Z/C}+\Delta_Z^{>0}))\otimes 
\mathcal A^k\right)
$$ 
is generated by its global sections for some positive integer $N$. 
Then $$
g_*\mathcal O_Z(k(K_{Z/C}+
\Delta_Z^{>0}))\otimes \mathcal A^{k-1}
$$ 
is nef on $C$. 
\end{lem}
\begin{proof}[Proof of Lemma \ref{ff-lem4.5}]
By the definition of $\mathcal F$ and \eqref{ff-eq4.1}, we have 
$$\psi_*\mathcal O_Z
(k(K_{Z/C}+\Delta^{>0}_Z))\simeq \psi_*\mathcal F. 
$$ 
This implies that 
$$g_*\mathcal O_Z
(k(K_{Z/C}+\Delta^{>0}_Z))\simeq g_*\mathcal F. 
$$ 
Therefore, $S^N(g_*\mathcal F\otimes \mathcal A^k)$ is generated by its global 
sections by 
assumption. 
Hence $|(\mathcal F\otimes g^*\mathcal A^k)^N|$ 
is a free linear system on $Z$. 
Note that $\mathcal F$ is $g$-generated. We set 
$$
\mathcal L=
\mathcal O_Z\left(K_{Z/C}
+\Delta^{=1}_{Z}+k\{\Delta^{>0}_Z\}\right)\otimes g^*\mathcal A. 
$$ 
Then we have 
$$
\mathcal L^k=
\mathcal O_Z
\left(E+(k-1)k\{\Delta^{>0}_Z\}\right)\otimes \mathcal F
\otimes g^*\mathcal A^k. 
$$ 
Let  $H$ be a 
general member of the free linear system $|(\mathcal F\otimes g^*\mathcal A^k)^N|$. 
Then we obtain 
$$
\mathcal L^{kN}=\mathcal O_Z
\left(H+NE+N(k-1)k\{\Delta^{>0}_Z\}\right). 
$$ 
We take a $(kN)$-fold cyclic cover $p:\widetilde Z\to Z$ associated to 
$$
\mathcal L^{kN}=\mathcal O_Z
\left(H+NE+N(k-1)k\{\Delta^{>0}_Z\}\right). 
$$ 
Note that $(\widetilde Z, p^*\Delta^{=1}_Z)$ is a semi-log-canonical 
pair (see \cite[Theorem 24]
{kollar-s}). More explicitly, $\widetilde Z$ can be written as follows: 
$$
\widetilde Z=\Spec _Z 
\bigoplus_{i=0}^{kN-1} \left(\mathcal L^{(i)}\right)^{-1}, 
$$ 
where 
$$
\left(\mathcal L^{(i)}\right)^{-1}=\mathcal L^{-i}\otimes \mathcal O_Z\left(
\left\lfloor \frac{i}{k}
\left(E+(k-1)k\{\Delta^{>0}_Z\}\right)\right\rfloor\right). 
$$ 
For the details of cyclic covers, see, for example, \cite[\S 3]{ev}, although 
\cite[\S 3]{ev} only treats the case where $Z$ is smooth. Since 
\begin{equation*}
p_*\omega_{\widetilde Z}\simeq \mathcal {H}om_{\mathcal O_Z}
\left (p_*\mathcal O_{\widetilde Z}, \omega_Z\right)\simeq 
\mathcal {H}om _{\mathcal O_Z} \left(\bigoplus _{i=0}^{kN-1} \left(\mathcal L^{(i)}\right)^{-1}, 
\omega_Z\right)\simeq \bigoplus _{i=0}^{kN-1} \omega_Z\otimes \mathcal L^{(i)}, 
\end{equation*} 
$\omega_Z\otimes \mathcal L^{(k-1)}$ is a direct summand of 
$p_*\omega_{\widetilde Z}$, 
where 
$$
\mathcal L^{(k-1)}=\mathcal L^{k-1}\otimes \mathcal O_Z
\left(-\left\lfloor \frac{k-1}{k}
\left(E+(k-1)k\{\Delta^{>0}_Z\}\right)\right\rfloor\right). 
$$ 
We note that 
\begin{align*}
&\mathcal O_Z\left(K_{Z/C}+\Delta ^{=1}_Z\right)\otimes 
\mathcal L^{(k-1)} \\&= 
\mathcal O_Z \left(k\left(K_{Z/C}+\Delta^{>0}_Z\right)
-\left \lfloor \frac{k-1}{k} E\right\rfloor\right) \otimes g^*\mathcal A^{k-1} \\ 
& \subset \mathcal O_Z \left(k\left(K_{Z/ C} 
+\Delta^{>0}_Z\right)\right)\otimes g^*\mathcal A^{k-1}
\end{align*} 
and that $E$ is the relative base locus of $\mathcal O_Z 
(k(K_{Z/C}+\Delta^{>0}_Z))$. 
Therefore, we have a natural isomorphism 
$$
g_*\left(\mathcal O_Z
(K_{Z/C}+\Delta^{=1}_Z)\otimes \mathcal L^{(k-1)}\right) 
\simeq g_*\mathcal O_Z
(k(K_{Z/C}+\Delta^{>0}_Z))\otimes \mathcal A^{k-1}
$$ 
since $$\left\lfloor \frac{k-1}{k}E\right\rfloor\leq E$$ (see 
the proof of \cite[Lemma 8.2]{fujino-weak}). 
We note that 
$\mathcal O_Z(K_Z+
\Delta^{=1}_Z)\otimes \mathcal L^{(k-1)}$ 
is a direct summand of 
$p_*\mathcal O_{\widetilde Z}(K_{\widetilde Z}+p^*\Delta ^{=1}_Z)$.  
By a special case of Theorem \ref{ff-thm1.10}, 
we obtain that $(g\circ p)_*\omega_{\widetilde Z}(p^*\Delta^{=1}_Z)$ is nef. 
Therefore, the direct summand 
$g_*\left(\mathcal O_Z(k(K_{Z/C}+\Delta^{=1}_Z)\otimes 
\mathcal L^{(k-1)}\right)$ is also a nef locally free sheaf on $C$. 
This means that $g_*\mathcal O_Z(k(K_{Z/C}+\Delta^{>0}_{Z}))
\otimes \mathcal A^{k-1}$ is a nef locally free sheaf on $C$. 
\end{proof}
By the definition of $r$, 
$g_*\mathcal O_Z(k(K_{Z/C}+\Delta_Z^{>0}))
\otimes \mathcal H^{rk-1}$ is nef. 
Therefore, $$
S^N\left(g_*\mathcal O_Z(k(K_{Z/C}
+\Delta_Z^{>0}))\otimes \mathcal H^{rk}\right)
$$ 
is generated by its global sections for some positive integer $N$. 
Then, by Lemma \ref{ff-lem4.5}, 
we obtain that  
$$
g_*\mathcal O_Z(k(K_{Z/C}+\Delta_Z^{>0}))\otimes 
\mathcal H^{r(k-1)}
$$ 
is nef. 
This is only possible if $(r-1)k-1<r(k-1)$. It is 
equivalent to $r\leq k$. 
Therefore, $$g_*\mathcal O_Z(k(K_{Z/C}+\Delta_Z^{>0}))
\otimes \mathcal H^{k^2-1}$$ is nef. 
The same holds true if we take any base change by $\pi:C'\to C$ 
such that $\pi$ is a finite morphism 
from a smooth projective curve and ramifies only over general 
points of $C$. 
Therefore, $$g_*\mathcal O_Z(k(K_{Z/C}+
\Delta_Z^{>0}))$$ is nef (see Lemma \ref{ff-lem2.3}).   

Hence we obtain that $f_*\mathcal O_X(k(K_{X/C}+\Delta))$ is a nef locally free 
sheaf on $C$. 
Since $\mathcal O_X(kl(K_X+\Delta))$ is $f$-generated,  by replacing $k$ with $kl$ in the 
above arguments, we 
obtain that $f_*\mathcal O_X(kl(K_{X/C}+\Delta))$ is nef for every positive integer $l$. 
\end{proof}

We close this section with comments on 
Koll\'ar's arguments in \cite[Section 4]{kollar} for 
the reader's convenience. 

\begin{say}[{Comments on Koll\'ar's arguments in \cite[Section 4]{kollar}}]
\label{ff-say4.6} 
In \cite[Section 4]{kollar}, 
Koll\'ar essentially claims Theorem \ref{ff-thm1.7} when
$\dim X=3$.
However, it is not so obvious to follow his arguments.
In the last part of \cite[4.14]{kollar}, he says: 
\begin{quote}
As in the proof of 4.13 the kernel of
$\delta$ is a direct summand and thus semipositive.
\end{quote}
In \cite[4.14]{kollar}, $E$ is not always smooth.
Therefore, it is not clear what kind of variations of 
Hodge structure should be considered.
The map
$$
(f\circ g)_*\omega_{E/C}\overset{\delta}{\longrightarrow} R^1(f\circ g)_*\omega_{X/C}
$$
in \cite[4.14]{kollar} is different from the map
$$
\delta': (f\circ g)_*\omega_{D'/C}\to R^1(f\circ g)_*\omega_{Z'/C}
$$
in the proof of \cite[4.13.~Lemma]{kollar} from the Hodge theoretic viewpoint.
Note that $D'$ and $Z'$ are smooth by construction.
In general, $E$ and $X$ are {\em{singular}} in \cite[4.14]{kollar}.
Koll\'ar informed the author that
the nefness of $(f\circ g)_*\omega_{X/C}(E)$ can be checked
with the aid of the classification of semi-log-canonical
surface singularities.
Note that his arguments only work for the case where 
the fibers are surfaces.
Anyway, we do not pursue them here because
the nefness of $(f\circ g)_*\omega_{X/C}(E)$ is a
special case of Theorem \ref{ff-thm1.10} when $f\circ g$ is projective. 
\end{say}

\section{Proof of corollaries}\label{ff-sec5}
In this final section, we prove the corollaries in Section \ref{ff-sec1}. 

\begin{proof}[Proof of Corollary \ref{ff-cor1.4}] 
The moduli functor $\mathcal M^{stable}_H$ is bounded by 
\cite{hmx}. Therefore, 
we have the coarse moduli space $M^{stable}_H$ of $\mathcal M^{stable}_H$ 
which is a complete algebraic space. 
Note that the moduli functor 
$\mathcal M^{stable}_H$ satisfies the valuative criterion of separatedness 
and the valuative criterion of properness. 
Since the moduli functor $\mathcal M^{stable}_H$ is semipositive in the 
sense of Koll\'ar by Theorems \ref{ff-thm1.2} and 
\ref{ff-thm1.7}, we see that $M^{stable}_H$ is a projective 
algebraic scheme (see \cite{kollar}). 
\end{proof}

Corollaries \ref{ff-cor1.5} and \ref{ff-cor1.6} are easy consequences of 
Corollary \ref{ff-cor1.4}. 

\begin{proof}[Proof of Corollary \ref{ff-cor1.5}] 
The moduli functor $\mathcal M_H$ is an open subfunctor 
of $\mathcal M^{stable}_H$ because the smoothness is an open condition. 
Therefore Corollary \ref{ff-cor1.5} follows from Corollary \ref{ff-cor1.4}. 
\end{proof}

\begin{proof}[Proof of Corollary \ref{ff-cor1.6}]
Note that any small deformations of canonical singularities are canonical 
(see \cite[Main Theorem]{kawamata}). 
Therefore the moduli 
functor $\mathcal M^{can}_H$ is an open subfunctor of 
$\mathcal M^{stable}_H$. 
Hence, Corollary \ref{ff-cor1.6} follows from Corollary \ref{ff-cor1.4}.  
\end{proof}


\end{document}